\documentclass[a4paper,11pt]{amsart}
\usepackage{amsmath,amsthm,amssymb}
\usepackage{comment}
\newcommand{\mychoice}[3]{#1
}
\mychoice{
\newcommand{\plabel}[1]{ \label{#1}}
\newcommand{\gbibitem}[1]{ \bibitem{#1}}

\specialcomment{commentx}{}{}\excludecomment{commentx}
\specialcomment{commenty}{}{}\excludecomment{commenty}
\newcommand{\rechoicecomm}[1]{}
}{
\usepackage{xcolor}
\newcommand{\plabel}[1]{ \label{#1}\rlap{\smash{${}^{^{[#1]}}$}}}
\newcommand{\gbibitem}[1]{ \bibitem{#1}\rlap{\smash{${}^{^{[#1]}}$}}}

\newenvironment{commenty}{\color{blue} }{\color{black} }
\newcommand{\rechoicecomm}[1]{#1}
}{
\usepackage{xcolor}
\newcommand{\plabel}[1]{ \label{#1}}
\newcommand{\gbibitem}[1]{ \bibitem{#1}}

\specialcomment{commentx}{}{}\excludecomment{commentx}

\newcommand{\rechoicecomm}[1]{}
}

\DeclareMathOperator{\diam}{diam}

\theoremstyle{definition}
\newtheorem{example}{Example}
\theoremstyle{plain}
\newtheorem{theorem}[example]{Theorem}

\newcommand{\eqedexer}{
\renewcommand{\qedsymbol}{$\diamondsuit$}
\pushQED{\qed}
\qedhere
\popQED
\renewcommand{\qedsymbol}{$\Box$}
}

\newcommand{\marginextend}[1]{ \addtolength{\oddsidemargin}{-#1}  \addtolength{\evensidemargin}{-#1}
  \addtolength{\textwidth}{#1}\addtolength{\textwidth}{#1}}
\newcommand{\updownextend}[1]{ \addtolength{\topmargin}{-#1}  \addtolength{\textheight}{#1}
\addtolength{\textheight}{#1}}
\marginextend{1cm}
\updownextend{0cm}
\allowdisplaybreaks[4]
\title{On averaged self-distances in finite dimensional Banach spaces}
\author{Gyula Lakos}
\address{Alfréd Rényi Institute of Mathematics, Reáltanoda utca 13-15, Budapest, H--1053, Hungary}
\email{lakos@renyi.hu}
\begin{document}
\begin{abstract}
Assume that $\mathfrak A$ is a real Banach space of finite dimension $n\geq2$.
Consider any Borel probability measure $\nu$ supported on the unit ball  $K$ of $\mathfrak A$.
We show that
\[\Delta(\nu)=\int_{x \in K}\int_{ y\in K}|x-y|_{\mathfrak A} \,\,\,\nu(x)\,\nu(y)\leq 2(1-2^{-n}f(n)),\]
where $f:\mathbb N\setminus \{0,1\}\rightarrow (0,1]$ is a concrete universal function such that $f(n)\sim \frac{2}{\mathrm e n^2\log n}$.
It is hoped that in the estimate`$f(n)$' can be replaced by `$1$'.
\end{abstract}
\maketitle
Assume that $\mathfrak A$ is a real Banach space of finite dimension $n$.
Consider any Borel probability measure $\nu$ supported on the unit ball  $K$ of $\mathfrak A$.
The averaged self-distance of $\nu$ is
\[\Delta(\nu)=\int_{x \in K}\int_{ y\in K}|x-y|_{\mathfrak A} \,\nu(x)\,\nu(y),\]
where $|\cdot|_{\mathfrak A}$ denotes the norm on $\mathfrak A$.
The trivial estimate is, of course, $\Delta(\nu)\leq2$.
This, in some cases, can be approached relatively well:
\begin{example}
Assume that $\mathfrak A$  is $\mathbb R^n$ endowed by the maximum norm $|(x_1,\ldots, x_n)|_\infty=\max(|x_1|,\ldots, |x_n|)$;
 and $\nu$ is distributed uniformly on the vertices $(\pm1,\ldots,\pm1)$ of the unit ball.
Then
\[\Delta(\nu)=2(1-2^{-n}).\eqedexer\]
\end{example}
In what follows we will be interested in estimates for $\Delta(\nu)$ smaller than $2$, depending on the dimension $n$,
 with emphasis on the asymptotical behaviour.
(Regarding the uniform distribution,
Arias-de-Reyna, Ball, Villa \cite{ABV}
investigates a similar question,  but from the other direction.)
For a warm-up, let us start with the case $n=2$.
\begin{theorem} For $n=2$,
 \[\Delta(\nu)\leq\frac{48}{31}+\frac{6}{31}\sqrt2=1.8221\ldots <\frac32+\frac14\sqrt2=1.8535\ldots. \]
\begin{proof}
By a theorem of Lassak \cite{L}, $K$ can be covered by four of its homothetical copies of ratio $\frac{\sqrt2}{2}$.
Let $K_1,\ldots,K_4$ be such copies.
Set $L_1=K\cap K_1$, $L_2=(K\cap K_2)\setminus K_1$, $L_3=(K\cap K_3)\setminus (K_1\cup K_2)$,
 $L_4=(K\cap K_4)\setminus(K_1\cup K_2\cup K_3)$.
Let $v_i=\nu(L_i)$.
Clearly, $v_1+\ldots+v_4=1$.

For $x,y\in L_i$, we can estimate $|x-y|_{\mathfrak A}$ by $\diam_{\mathfrak A} K_i=\sqrt2$,
 and it can be estimated by $2$ otherwise.
Therefore,
\[\Delta(\nu)\leq 2(v_1+\ldots+v_4)^2-(2-\sqrt2)( (v_1)^2+\ldots+(v_4)^2 ).\]
Using the inequality between the arithmetical and square means, we obtain
\[\leq2(v_1+\ldots+v_4)^2-(2-\sqrt2) \frac14 (v_1+\ldots+v_4)^2=2-(2-\sqrt2) \frac14=\frac32+\frac14\sqrt2.\]

One can do slightly better.
Let $\varkappa_2$ be the least upper bound achievable for $\Delta(\nu)$ in general.
Then this bound can be applied for $\nu|_{L_i}$ supported on $K_i$ except it must be taken into consideration that
 the restricted measures are not probability measures and that $K_i$ is a scaled version of $K$.
In this manner, we obtain
\[\Delta(\nu)\leq 2(v_1+\ldots+v_4)^2- \left(2-\varkappa_2\frac{\sqrt2}2\right) ( (v_1)^2+\ldots+(v_4)^2 )
\leq2- \left(2-\varkappa_2\frac{\sqrt2}2\right) \frac14.\]
This also applies if $\Delta(\nu)$ approaches $\varkappa_2$; therefore
$\varkappa_2\leq2-\left(2-\varkappa_2\frac{\sqrt2}2\right) \frac14,$
which leads to  $\varkappa_2\leq\frac{48}{31}+\frac{6}{31}\sqrt2$.
(But even so, this estimate is very crude.)
\end{proof}
\end{theorem}
The higher dimensional case can be treated similarly.
Let us recall certain estimates due to Rogers.
Let us consider any centrally symmetric compact convex body $H$ in the $n$-dimensional
 space, $\mathbb R^n$.
Let $\Theta(H)$ be the infimum of the covering density of $\mathbb R^n$ by translates
 of $H$; and let $\Theta_L(H)$ be the same with respect to lattice coverings.
Then, by Rogers \cite{R}, for $n\geq3$,
\begin{align}
\Theta(H)\leq \Theta_L(H)\leq&\min_{0<\eta<1/n}(1+\eta)^n(1+n\log(1/\eta))\plabel{eq:r1}\\
<& \left(1+\frac1{n\log n}\right)^n(1+n\log(n\log n))\plabel{eq:r2}\\
<& n\log n +n\log\log n+2n+1 \plabel{eq:r3}\\
<& n\log n +n\log\log n+5n.\plabel{eq:r4}
\end{align}
Here  \eqref{eq:r1} is the proper result of \cite{R} (which can be expressed explicitly using the Lambert $W_{-1}$ function);
 \eqref{eq:r2} reflects the choice $\eta=\frac1{n\log n}$;
and while, for example,  \eqref{eq:r3} is still true (but we do no want to spend time on it),
 \eqref{eq:r4}  is the generally quoted estimate; see G. Fejes Tóth \cite{F} for more on this.
The main point is that there is a dimension-dependent but otherwise universal quantity $\Theta_n$,
which is nearly linear in $n$, estimating the minimal translative covering densities from above.
Furthermore, as it is explained in Rogers, Zong \cite{RZ}, if $0<r<1$, then $H$ can be covered by at
most $(1+r^{-1})^n\Theta(H)\leq (1+r^{-1})^n\Theta_n$
translates of $rH$ (i. e. by homothetical copies of $H$ with ratio $r$).
For a recent review on these topics, see Naszódi \cite{N}.
\begin{theorem} For $n\geq3$,
 \[\Delta(\nu)<2\left(1-2^{-n}\frac{2\left(1-\frac2n\right)}{\mathrm e n \Theta_n} \right). \]
\begin{proof}
Let us apply the argument of  Rogers, Zong \cite{RZ} with $H=  K$ in $\mathfrak A\simeq\mathbb R^n$.
Then $K$ can be covered by at most $s= (1+r^{-1})^n $ many copies of $rK$, say $K_1,\ldots, K_{\lfloor s\rfloor}$.
Let $L_1=K\cap K_1$, and $L_i=(K\cap K_i )\setminus  (K_1\cup \ldots \cup K_{i-1})$.
Set $v_i=\nu(L_i)$.
Then
\begin{multline*}
\Delta(\nu)\leq 2(v_1+\ldots+v_{\lfloor s\rfloor})^2- \left(2-2r\right) ( (v_1)^2+\ldots+(v_{\lfloor s\rfloor})^2 ),
\\
\leq2-2(1-r)\frac{1}{\lfloor s\rfloor}\leq 2-2(1-r)\frac{1}{ s }=2-2\frac{(1-r)(1+r^{-1})^{-n}}{\Theta_n} .
\end{multline*}
Now $r$ can be optimized to $r=\frac12\left(\sqrt{n^2+6n+1}-n-1\right)$, but, for the sake of simplicity, we take $r=1-\frac2n$.
With this latter choice,
\[\Delta(\nu)\leq 2\left(1- {\frac {{2}^{1-n}}{n} \left( {\frac {n-2}{n-1}} \right) ^{ n}}\frac1{\Theta_n}\right)<
 2\left(1- {\frac {{2}^{1-n}}{n} \left( {\frac {1-\frac2n}{\mathrm e}} \right)  }\frac1{\Theta_n}\right).
\qedhere\]
\end{proof}
\end{theorem}
Altogether, this means that
$\Delta(\nu)\leq 2(1-2^{-n}f(n))$, where  $f:\mathbb N\setminus \{0,1\}\rightarrow (0,1]$ is an
explicit function with  $f(n)\sim \frac{2}{\mathrm e n^2\log n}$.

What one would expect is that $\Delta(\nu)\leq2(1-2^{-n})$ holds in general, with equality
if and only if $K$ is a parallelepiped and $\nu$ is distributed uniformly on its vertices.
That would be a continuous counterpart of Hadwiger's conjecture, cf. Hadwiger \cite{H},
Boltyanski,  Martini,  Soltan \cite{BMS}.

\textit{Acknowledgements.}
The author is extremely grateful to Márton Naszódi and Endre Makai, Jr. for helpful discussions and comments.

\end{document}